\documentclass[11pt]{amsart}
\usepackage{amssymb,latexsym,amsmath}
\input{xypic}
\xyoption{all}
\numberwithin{equation}{section}

\newcommand{\Q}{{\mathbb Q}}
\newcommand{\R}{{\mathbb R}}
\newcommand{\Z}{{\mathbb Z}}
\newcommand{\C}{{\mathbb C}}

\newcommand{\A}{{\mathbb A}}

\newcommand{\g}{{\mathfrak g}}

\newcommand{\M}{{\mathcal M}}

\newcommand{\GL}{{\rm GL}}
\newcommand{\SO}{{\rm SO}}

\newcommand\SGK{\mathcal{S}^G_{K_f}}
\newcommand\HGK{\mathcal{H}^G_{K_f}}

\newcommand\HG{\mathcal{H}_G^\place}
\newcommand\tM{\widetilde{\mathcal{M}}}
\def\fR{\mathfrak{R}}
\newcommand\place{\mathsf{S}}
\def\r{\mathfrak{r}}

\def\tM{\widetilde{\mathcal{M}}}

\def\aInd{{}^{\rm a}{\rm Ind}}

\numberwithin{equation}{section}
\newtheorem{theorem}[equation]{Theorem}

\newtheorem{lemma}[equation]{Lemma}

\begin{document}

\title[$L$-functions for orthogonal groups]{Special values of $L$-functions for orthogonal groups}
\date{\today}
\thanks{C.B. is partially supported by DST-INSPIRE Faculty scheme, award number [IFA- 11MA -05].}
\author{\bf Chandrasheel Bhagwat \ \ \& \ \ A. Raghuram}
\address{Indian Institute of Science Education and Research, Dr.\,Homi Bhabha Road, Pashan, Pune 411008,  INDIA.}
\email{cbhagwat@iiserpune.ac.in, \ raghuram@iiserpune.ac.in}

\begin{abstract} 
This is an announcement of certain rationality results for the critical values of the 
degree-$2n$ $L$-functions attached to 
$\GL_1 \times \SO(n,n)$ over $\Q$ for an even positive integer $n$.  The proof follows from studying the rank-one Eisenstein cohomology for 
$\SO(n+1,n+1).$  

\end{abstract}

\maketitle

\section{\bf Introduction and statement of the main result}
To motivate the main result, let's recall a well-known theorem of Shimura \cite{shimura-mathann}. 

\begin{theorem}[Shimura]
Let $f = \sum a_nq^n \in S_k(N,\chi)$ and $g = \sum b_nq^n \in S_l(N,\psi)$ be primitive modular forms of weights 
$k$ and $l$, with nebentypus characters $\chi$ and $\psi$ for $\Gamma_0(N).$
Let $\Q(f,g)$ be the number field obtained by adjoining the Fourier coefficients $\{a_n\}$ and $\{b_n\}$ to $\Q.$ 
Assume that $k > l.$ Let 
$$
D_N(s, f,g) \ := \ L_N(2s+2-k-l, \chi\psi) \sum_{n=1}^\infty \frac{a_nb_n}{n^s} 
$$
be the degree $4$ Rankin--Selberg L-function attached to the pair $(f,g)$. Then, for any integer $m$ with 
$l \leq m < k$ we have: 
$$
D_N(m,f,g) \ \approx \ (2\pi i)^{l+1-2m} \, \g(\psi) \, u^+(f)u^-(f),  
$$
where $\approx$ means up to an element of $\Q(f,g)$, $u^\pm(f)$ are the two periods attached to $f$ by Shimura, and 
$\g(\psi)$ is the Gauss sum of $\psi.$ Furthermore, the ratio of the $L$-value in the left hand side by the right hand side is equivariant under ${\rm Gal}(\overline{\Q}/\Q).$
\end{theorem}

\bigskip

The integers $l \leq m < k$ are all the critical points for $D_N(s, f,g).$ (There are no critical points if $l = k.$) 
Suppose $k \geq l+2$, and we look at two successive critical values then the only change in the right hand side is 
$(2\pi i)^{-2}$ which may be seen to be exactly accounted for 
by the $\Gamma$-factors at infinity. Suppose $L(s, f \times g)$ denotes the completed degree-$4$ $L$-function attached to $(f,g)$, normalized in a classical way as in the theorem above, then we deduce: 
\begin{equation}
\label{eqn:shimura}
L(l, f \times g) \ \approx \ L(l+1, f \times g) \ \approx \ \cdots \ \approx \ L(k-1, f \times g). 
\end{equation}

\medskip

The above result is a statement for $L$-functions for $\GL_2 \times \GL_2$ over $\Q.$ 
Later Shimura generalized this to Hilbert modular forms \cite{shimura-duke}, i.e., for $\GL_2 \times \GL_2$ over a totally real field $F.$ Note that $(\GL_2 \times \GL_2)/\Delta\GL_1 \simeq {\rm GSO}(2,2),$ i.e., Shimura's result may be construed as a theorem for $L$-functions for orthogonal groups in four variables.

{\it The main aim of this article is to announce that we can generalize the result in (\ref{eqn:shimura}) to $L$-functions for 
$\GL_1 \times \SO(n,n)$ over a totally real field $F$, and when $n = 2r \geq 2.$} For simplicity of exposition, we will work over $F = \Q.$

\bigskip

\begin{theorem}
\label{thm:main}
Let $n = 2r \geq 2$ be an even positive integer. 
Consider $\SO(n,n)/\Q$ defined so that the subgroup of all upper-triangular matrices is a Borel subgroup. 
Let $\mu$ be a dominant integral weight written as $\mu = (\mu_1 \geq \mu_2 \geq \dots \geq \mu_{n-1} \geq |\mu_n|),$ with $\mu_j \in \Z.$ 
Let $\sigma$ be a cuspidal automorphic representation of $\SO(n,n)/\Q.$ Assume: 

\begin{enumerate}
\item the Arthur parameter $\Psi_\sigma$ is cuspidal on $\GL_{2n}/\Q$; 
\item $\sigma$ is globally generic; 
\item $\sigma_\infty|_{\SO(n,n)(\R)^{0}}$ is a discrete series representation with Harish-Chandra parameter 
$\mu + \rho_n.$
\end{enumerate}
Let ${}^\circ\chi$ be a finite order character of $\Q^\times\backslash \A^\times.$ 
Then the critical set for the degree-$2n$ completed $L$-function $L(s, {}^\circ\chi \times \sigma)$ is the finite set of contiguous integers 
$$
\{ 1-|\mu_n|, \, 2-|\mu_n|, \, \dots, \,  |\mu_n|\}.
$$
Assume also that $|\mu_n| \geq 1$, so that the critical set is nonempty; and in this case there are at least two critical points. We have
$$
L(1-|\mu_n|, {}^\circ\chi \times \sigma) \ \approx \ L(2-|\mu_n|, {}^\circ\chi \times \sigma) \ \approx \ \cdots \ \approx \ 
L(|\mu_n|, {}^\circ\chi \times \sigma), 
$$
where $\approx$ means up to an element of a number field $\Q({}^\circ\chi, \sigma)$, and furthermore, all the successive ratios are equivariant under ${\rm Gal}(\overline{\Q}/\Q).$
\end{theorem}

\bigskip

\section{\bf The combinatorial lemma and a restatement of the main theorem} 

The strategy of proof follows the paradigm in Harder-Raghuram~\cite{harder-raghuram-CR} \cite{harder-raghuram}. In 
our situation, this  involves studying the rank-one Eisenstein cohomology of $G = \SO(n+1,n+1)$, especially the contribution coming from a parabolic subgroup $P$ with Levi quotient $M_P = \GL_1 \times \SO(n,n).$ 
As in {\it loc.\,cit.}\,certain Weyl group combinatorics play an important role--essentially saying that a particular context involving the cohomology of arithmetic groups is viable exactly when the intervening $L$-values are critical.

\begin{lemma}
Let $\mu = (\mu_1 \geq \mu_2 \geq \dots \geq \mu_{n-1} \geq |\mu_n|)$ be a dominant integral weight, and  
$\sigma$ be a cuspidal automorphic representation for $\SO(n,n)/\Q$ as in Thm.\,\ref{thm:main}. 
Let $d \in \Z$ and put $\chi = |\ |^{-d} \otimes {}^\circ\chi$ where ${}^\circ\chi$ is a finite-order character.
Let $G = \SO(n+1,n+1)$ and $P$ the maximal parabolic subgroup obtained by deleting the `first' simple root, in which case the Levi decomposition $P = M_P N_P$ looks like: $M_P = \GL_1 \times \SO(n,n)$ and ${\rm dim}(N_P) = 2n.$ 
The following are equivalent: 
\medskip
\begin{enumerate}
\item $-n$ and $1-n$ are critical for the completed degree-$2n$ $L$-function $L(s, \chi \times \sigma).$ 
\medskip
\item $1- |\mu_n| \ \leq \ n+d \ \leq \ |\mu_n| -1.$ 
\medskip
\item There is a unique $w \in W^P$ (here $W^P$ is the set of Kostant representatives for $P$; we have 
$W_G = W_{M_P}W^P$) such that $w^{-1} \cdot (d \times \mu)$ is dominant for $G$ and $l(w) = {\rm dim}(N_P)/2.$ 
\end{enumerate}
\end{lemma}

As $d$ runs through the range prescribed by (2), the ratio of critical values 
$$
\frac{L(-n, \chi \times \sigma)}{L(1-n, \chi \times \sigma)}
$$
(where the criticality is assured by (1)) runs through the set of all successive ratios of critical values 
$$
\left\{
\frac{L(1-|\mu_n|, {}^\circ\chi \times \sigma)}{L(2-|\mu_n|, {}^\circ\chi \times \sigma)}, \ \dots, \ 
\frac{L(|\mu_n|-1, {}^\circ\chi \times \sigma)}{L(|\mu_n|, {}^\circ\chi \times \sigma)}
\right\}. 
$$

\medskip

This says that when the method of Eisenstein cohomology is invoked for rationality results, 
then we get a result  for ratios of all possible successive critical values, no more and no less! It suffices now to prove the following

\begin{theorem}
\label{thm:main-restated}
Let the notations on $\chi$ and $\sigma$ 
be as in the lemma above, and assume that the conditions on $d$ are satisfied. Then the ratio of 
$L$-values 
$$
\frac{L(-n, \chi \times \sigma)}{L(1-n, \chi \times \sigma)} 
$$
is algebraic and is ${\rm Gal}(\overline{\Q}/\Q)$-equivariant.
\end{theorem}

\bigskip

\section{\bf Comments on the consequences of various hypotheses of the main theorem}

\subsection{A discrete series representation as the local representation at infinity}
\label{sec:discrete-galois} 
This is the simplest kind of representation with nontrivial relative Lie algebra cohomology; in fact, it has nonzero cohomology only in the middle degree. Furthermore, this implies that the finite part $\sigma_f$ contributes to the cohomology of a locally symmetric space of $\SO(n,n)$ with coefficients in 
the local system attached to $\mu.$ Using arguments as in Gan-Raghuram~\cite{gan-raghuram}, we show that 
$\sigma_f$ is defined over a number field $\Q(\sigma)$ and there is a ${\rm Gal}(\bar\Q/\Q)$-action on the set of cuspidal 
representations that satisfy the hypotheses (1), (2) and (3).  In the statement of the theorem above, 
$\Q({}^\circ\chi, \sigma)$ is the field generated by the values of ${}^\circ\chi$ and $\Q(\sigma)$.

\subsection{The transfer $\Psi_\sigma$ is cuspidal on $\GL_{2n}/\Q$}
This is needed for two reasons: (1) We do not want the $L$-function 
$L(s, {}^\circ\chi\times\sigma)$ to break up into smaller $L$-functions; although, even if it did, with an inductive argument, at least in the case when $\Psi_\sigma$ is tempered, we would very likely still have the main theorem. 
(2) The second reason is far more serious and very delicate. We need to prove a `Manin-Drinfeld' principle: that there is a 
Hecke-projection from the total boundary cohomology (of the Borel-Serre boundary) 
to the isotypic component of the representation of $G$ induced from $\chi \otimes \sigma$ of $M_P.$ 
See Sect.\,\ref{sec:proof} below. 
For this to work, we have to exclude the possibility of $\sigma$ being, for example, a CAP representation (which also gets guaranteed by the next hypothesis).

\subsection{$\sigma$ is globally generic} 
This hypothesis plays several roles: it is used in proving the existence of Galois action mentioned in 
Sect.\,\ref{sec:discrete-galois} above.  
Shahidi's results \cite{shahidi-duke85} on local constants (see Sect.\,\ref{sec:proof} below)  
need genericity of the representation at infinity.

\subsection{Compatibility with Deligne's conjecture} 
The above theorem would be compatible with Deligne's conjecture \cite{deligne} on the critical values of motivic $L$-function, if we have the following period relation: 
{\it Let $M$ be a pure regular motive of rank-$2n$ over $\Q$ with coefficients in a number field $E.$ Suppose $M$ is of 
orthogonal type (i.e., there is a map ${\rm Sym^2}(M) \to \Q(-{\sf w})$, where ${\sf w}$ is the purity weight of $M$), then
Deligne's periods $c^\pm(M)$ are related as:}
$$
c^+(M) \ = \ c^-(M), \ \mbox{\it as elements of $(E \otimes \C)^\times/E^\times.$}
$$
This is known if $M$ is a tensor product of two rank-two motives; see Blasius~\cite[2.3]{blasius}.

\subsection{Langlands transfer and special values} 
It is important to prove this theorem at the level of $L$-functions for $\GL_1 \times \SO(n,n)$, and {\it not} as 
$L$-functions for $\GL_1 \times \GL_{2n}$ after transferring. We would see this subtle point already in the context of Shimura's theorem, because (i) the Langlands transfer $f \boxtimes g$, which is a cuspidal representation of $\GL_4$ does not see the Petersson norm $\langle f, f \rangle$ of only one of the constituents; and (ii) for an 
$L$-function $L(s, \pi)$ with $\pi$ cuspidal on $\GL_4/\Q$, successive $L$-values would see $c^+(\pi)$ and $c^-(\pi)$, and in the automorphic world, it is not (yet) known that if $\pi$ came via transfer from $\GL_2 \times \GL_2$ then 
$c^+(\pi) \approx c^-(\pi).$  In a similar vein, one may ask if the main result of \cite{harder-raghuram} applied to 
$\GL_1 \times \GL_{2n}$ implies the main result of this paper; this would be so if we could prove that the relative periods,
denoted $\Omega^\varepsilon$ therein, for the representation $\Psi_\sigma$ of $\GL_{2n}$ are trivial--at this moment we have no idea how one might prove such a period relation--hence our insistence on working intrinsically  in the context of orthogonal groups.

\subsection{Further generalizations} 
All this {\it should} work for $L$-functions for $\GL_1 \times {\rm GSpin}(2n)$ over a totally real field $F.$ We say {\it should} because of the hypothesis ``the Arthur parameter $\Psi_\sigma$ being cuspidal." We may appeal to the work of Asgari and Shahidi \cite{asgari-shahidi-duke} since we only want the case of {\it generic transfer} 
from ${\rm GSpin}(2n)$ to ${\rm GL}_{2n}.$ However, as we see below, this hypothesis is also needed for the Manin-Drinfeld principle for boundary cohomology, and for this we will need Arthur's work \cite{arthur}. Using the results in 
the recent thesis of  Bin Xu \cite{xu}, it might be possible to generalize our results to $\GL_1 \times {\rm GO}(n,n).$

\bigskip
\section{\bf An adumbration of the proof of Theorem~\ref{thm:main-restated}}
\label{sec:proof}

The basic idea, following \cite{harder-raghuram-CR} and \cite{harder-raghuram}, is to give a cohomological interpretation to the constant term theorem of Langlands, by studying the rank-one Eisenstein cohomology of $\SO(n+1,n+1)$. Let the notations be as in the combinatorial lemma above. A consequence of this lemma is that the representation algebraically (un-normalized) and parabolically induced from 
$\chi_f \otimes \sigma_f $ appears in boundary cohomology: 
$$
\aInd_{P(\A_f)}^{G(\A_f)}(\chi_f \otimes \sigma_f)^{K_f} \ \hookrightarrow 
H^{q_0}(\partial_P\SGK, \tM_{\lambda,E}), 
$$
where $q_0 = \mbox{middle-dimension-of-symmetric-space-of-$M_P$} + {\rm dim}(N_P)/2;$  
$\lambda = w^{-1} \cdot (d + \mu);$ $K_f$ is a deep-enough open-compact subgroup of $G(\A_f)$; 
$\partial_P$ denotes the part corresponding to $P$ of the Borel-Serre boundary  of the locally symmetric space 
$\SGK$ for $G$ with level structure $K_f$; 
$ \tM_{\lambda,E}$ is the sheaf corresponding to the finite-dimensional representation $\M_{\lambda,E}$ of the 
algebraic group $G \times E.$ (The reader is referred to \cite[Sect.\,1]{harder-raghuram-CR} for a quick primer on these cohomology groups and for the fundamental long exact sequence that comes out of the Borel-Serre compactification.)
The field $E$ is taken to be a large enough Galois extension of $\Q$, for example, take $E$ to contain $\Q(\chi, \sigma)$. 
To relate to the theory of automorphic forms, 
we can pass to $\C$ via an embedding $\iota : E \to \C.$ 
The induced representations and the cohomology groups are all modules for a Hecke-algebra $\HGK$, and in what follows below, we restrict our attention to a commutative sub-algebra $\HG$ ignoring a finite set $\place$ of all ramified places. 
Next, one observes that the standard intertwining operator $T_{\rm st}$, at the point of evaluation $s = - n$ goes as: 
$$
T_{\rm st} \ : \ \aInd_{P(\A_f)}^{G(\A_f)}(\chi_f \otimes \sigma_f) \ \longrightarrow \ 
 \aInd_{P(\A_f)}^{G(\A_f)}(\chi_f^{-1}(2n) \otimes {}^\kappa\!\sigma_f), 
$$
where $(2n)$ denotes a Tate-twist, and $\kappa$ is an element of ${\rm O}(n,n)$ but outside $\SO(n,n)$. 
Certain combinatorial details about Kostant representatives allow us to observe that the induced representation 
in the target also appears in boundary cohomology as: 
$$
\aInd_{P(\A_f)}^{G(\A_f)}(\chi^{-1}_f(2n) \otimes {}^\kappa\!\sigma_f)^{K_f} \ \hookrightarrow 
H^{q_0}(\partial_P\SGK, \tM_{\lambda,E}), 
$$
for the same degree $q_0$ and the same weight $\lambda.$ Let  
$$
I^\place_P(\chi_f, \sigma_f)^{K_f} \ := \ 
\aInd_{P(\A_f)}^{G(\A_f)}(\chi_f \otimes \sigma_f)^{K_f} \ \oplus \ 
 \aInd_{P(\A_f)}^{G(\A_f)}(\chi_f^{-1}(2n) \otimes {}^\kappa\!\sigma_f)^{K_f}. 
$$
The Manin-Drinfeld principle amounts to showing that we get a $\HG$-equivariant 
projection from 
boundary cohomology onto $I^\place_P(\chi_f, \sigma_f)^{K_f},$ 
and the target is isotypic, i.e., doesn't weakly intertwine 
with the quotient of the boundary cohomology by $I^\place_P(\chi_f, \sigma_f)^{K_f}.$ Denote this projection as: 
$$
\fR : H^{q_0}(\partial \SGK,\tM_{\lambda,E}) \ \longrightarrow \ 
I^\place_P(\chi_f, \sigma_f)^{K_f}. 
$$
If we denote the restriction map from global cohomology to the boundary cohomology as 
$\r^* : H^{q_0}( \SGK,\tM_{\lambda,E}) \to H^{q_0}(\partial \SGK,\tM_{\lambda,E})$, then the main technical result on Eisenstein cohomology involves the image of the composition $\fR \circ \r^*$:  
$$
H^{q_0}( \SGK,\tM_{\lambda,E})  \ \stackrel{\r^*}{\longrightarrow} \ 
H^{q_0}(\partial \SGK,\tM_{\lambda,E}) \ \stackrel{\fR}{\longrightarrow} \ 
I^\place_P(\chi_f, \sigma_f)^{K_f}. 
$$

\smallskip

For simplicity of explanation, let's pretend (and this could very well happen in some cases) that $I^\place_P(\chi_f, \sigma_f)^{K_f}$ is a two-dimensional $E$-vector space. Our main result on Eisenstein cohomology will then say that 
the image of $\fR \circ \r^*$ is a one-dimensional subspace of this two-dimensional ambient space. We then look at the slope of this line. Passing to a transcendental level via an $\iota : E \to \C$, and 
using the constant term theorem, one proves that the slope is in fact
$$
c_\infty(\chi_\infty,\sigma_\infty) \frac{L_f(-n, \chi \times \sigma)}{L_f(1-n, \chi \times \sigma)}, 
$$
where $c_\infty(\chi_\infty,\sigma_\infty)$ is a nonzero complex number depending only on the data at infinity, and 
$L_f(s, \chi \times \sigma)$ is the finite part of the $L$-function. This proves that above quantity lies in $\iota(E).$ 
Studying the behavior of the cohomology groups on varying $E$ then proves Galois-equivariance. 

\smallskip

Along the way, we need to address certain local problems. At the finite ramified places we prove that the local normalized intertwining operator is nonzero and preserves rationality using the results of Kim \cite{kim}, 
M\oe glin--Waldspurger \cite{moeglin-waldspurger} and Waldspurger \cite{waldspurger}. 
At the archimedean place, 
yet another consequence of the combinatorial lemma is that the representation 
$$
\aInd_{P(\R)}^{G(\R)}(\chi_\infty \otimes \sigma_\infty)
$$
is irreducible; this follows from the results of Speh-Vogan \cite{speh-vogan}. Using Shahidi's results \cite{shahidi-duke85} on local factors we then deduce that the standard intertwining operator is an isomorphism and induces a nonzero isomorphism in relative Lie algebra cohomology. But these cohomology groups at infinity are one-dimensional, and 
after fixing bases on either side we get a nonzero number $c_\infty(\chi_\infty,\sigma_\infty)$.  
We expect that a careful analysis, as in Harder~\cite{harder}, of the rationality properties of relative Lie algebra cohomology groups, should give us that $c_\infty(\chi_\infty,\sigma_\infty)$ is  the same as $L(-n, \chi_\infty \times \sigma_\infty)/L(1-n, \chi_\infty \times \sigma_\infty)$ up to a nonzero rational number, justifying our claim about a rationality result for completed $L$-values.  


\end{document}